\documentclass[preprint]{article}

\usepackage{amsmath}
\usepackage{amsthm}
\usepackage[all]{xy}

\newcommand{\ZZ}{\mathbf{Z}}

\theoremstyle{plain} 
\newtheorem{theorem}{Theorem}
\newtheorem{lemma}{Lemma}
\newtheorem{proposition}{Proposition}
\theoremstyle{definition}
\newtheorem{definition}{Definition}

\begin{document}

\title{The nilpotence degree of torsion elements in lambda--rings}

\author{F.J.-B.J.~Clauwens}

\maketitle

\section{Introduction.}

If $A$ is a $\lambda$-ring 
then for each prime number $p$
one has a map
$\theta^p\colon A\to A$ such that
\begin{itemize}
\item
The map $\psi^p\colon A\to A$ defined by
$\psi^p(a)=a^p-p\theta^p(a)$ is a ring homomorphism.
\item
The following formulas hold for all $a,b\in A$:
\begin{equation*}
\begin{split}
\theta^p(1)&=0\\
\theta^p(a+b)&=\theta^p(a)+\theta^p(b)
+\sum_{j=1}^{p-1}\frac{1}{p}\binom{p}{j}a^jb^{p-j}\\
\theta^p(ab)&=\theta^p(a)\psi^p(b)+a^p\theta^p(b)\\
\theta^p\psi^p(a)&=\psi^p\theta^p(a)
\end{split}
\end{equation*}
\end{itemize}
Moreover the maps $\psi^p$ and $\theta^r$ commute also for $r\not=p$.
Conversely if a commutative ring $A$ is equipped with maps
$\psi^p$ and $\theta^p$ satisfying the above identities
then $A$ has a unique structure of $\lambda$-ring such that the
$\psi^p$ are the associated Adams operations.
We refer to \cite{flog} for generalities about $\lambda$-rings.

In \cite{dress} it is shown that a torsion element in a $\lambda$-ring is nilpotent.
The purpose of this note is to give a sharp estimate of the nilpotence degree.

Let us call a ring $A$ equipped with one operation $\theta^p$
satisfying the above conditions a $\theta^p$-ring.
We will prove that if $A$ is a $\theta^p$-ring
and $a\in A$ satisfies $p^ea=0$ then $a^{p^e+p^{e-1}}=0$.
In order to show that this estimate is sharp we exhibit
a $\lambda$-ring $A$ and an element $a\in A$
such that $p^ea=0$ but $a^{p^e+p^{e-1}-1}\not=0$.

If $a$ is an element of a $\lambda$-ring $A$ and $na=0$
for some $n\in\ZZ$ then by the above theorem we have
\begin{equation*}
(np^{-e})^{p^e+p^{e-1}}a^{p^e+p^{e-1}}=0
\end{equation*}
where $p^e$ is the highest power of $p$ dividing $n$.
Since the greatest common divisor
of the numbers  $(np^{-e})^{p^e+p^{e-1}}$ is $1$ it follows that
\begin{equation*}
a^E=0,\text{ where }E=\max\{p^e+p^{e-1}\;;\;p^e\text{ divides }n\}
\end{equation*}

\section{The nilpotency estimate.}

In this section $A$ is a $\lambda$-ring at $p$, and $a\in A$
satisfies $p^ea=0$ for some $e>0$.
We will write $\psi^{p^e}$ for the $e$-fold iterate of $\psi^p$.
\begin{proposition}
For any $b\in A$ one has
$\theta^p(pb)=p^{p-1}b^p-\psi^p(b)$.
\end{proposition}

\begin{proof}
By definition of $\lambda$-ring at $p$ one has  $\theta^p(pb)=\theta^p(p)b^p+\psi^p(p)\theta^p(b)$.
Now substitute $\psi^p(p)=p$, $\theta^p(p)=p^{p-1}-1$
and $p\theta^p(b)=b^p-\psi^p(b)$.
\end{proof}

\begin{proposition}
\label{two}
For $0\leq k\leq e$ one has $p^{e-k}\psi^{p^k}(a)=0$.
\end{proposition}

\begin{proof}
By induction on $k$.
For $k=0$ the conclusion $p^ea=0$ is given.
Now assume the statement is true for $k=m$.
Then by the last proposition one has
\begin{equation*}
\begin{split}
0&=\theta^p(0)=\theta^p(p\cdot p^{e-m-1}\psi^{p^m}(a))\\
&=p^{p-1}(p^{e-m-1}\psi^{p^m}(a))^p-\psi^p(p^{e-m-1}\psi^{p^m}(a))\\
&=-\psi^p(p^{e-m-1}\psi^{p^m}(a))
=p^{e-m-1}\psi^{p^{m+1}}(a)\\
\end{split}
\end{equation*}
since $(p-1)+p(e-m-1)\geq e-m$ for $m<e$.
\end{proof}

\begin{proposition}
For $0\leq k\leq e-1$ one has
$a^{p^e+p^{e-1}}=
\psi^{p^k}(a^{p^{e-k}+p^{e-k-1}})$.
\end{proposition}

\begin{proof}
By induction on $k$.
For $k=0$ both sides are identical.
Assume that the statement is true for $k=m<e-1$.
Then one has
\begin{equation*}
\begin{split}
a^{p^e+p^{e-1}}
&=(\psi^{p^m}(a)^p)^{p^{e-m-2}(p+1)}\\
&=(\psi^p\psi^{p^m}(a)+p\theta^p\psi^{p^m}(a))^{p^{e-m-2}(p+1)}\\
&=\sum_i\binom{p^{e-m-2}(p+1)}{i}(p\theta^p\psi^{p^m}(a))^i
\psi^{p^{m+1}}(a)^{p^{e-m-2}(p+1)-i}
\end{split}
\end{equation*}
If $0<i<p^{e-m-2}(p+1)$ then the number of factors $p$
in the corresponding term is
\begin{equation*}
\begin{split}
&i+v_p\binom{p^{e-m-2}(p+1)}{i}
\geq i+v_p(p^{e-m-2}(p+1))-v_p(i)\\
&\geq i+(e+m-2)-(i-1)=e+m-1
\end{split}
\end{equation*}
and it also contains a factor $\psi^{p^{m+1}}(a)$.
Hence  by Proposition \ref{two} it vanishes.

If $i=p^{e-m-2}(p+1)$ then the term is a multiple of
$p^{e+1-m}\theta^p\psi^{p^m}(a)$, 
since $p^j(p+1)\geq j+3$ for all $j\geq 0$.
But this expression equals
\begin{equation*}
\begin{split}
&p^{e+1-m}\psi^{p^m}\theta^p(a)=
p^{e-m}\psi^{p^m}(a^p-\psi^p(a))=\\
&p^{e-m}\psi^{p^m}(a)^p-\psi^p(p^{e-m}\psi^{p^m}(a))=0-0=0
\end{split}
\end{equation*}
by Proposition \ref{two}.

Thus only the term with $i=0$ remains, which proves the statement 
for $k=m+1$.
\end{proof}

\begin{theorem}
$a^{p^e+p^{e-1}}=0$.
\end{theorem}

\begin{proof}
By the last Proposition we have
\begin{equation*}
\begin{split}
&a^{p^e+p^{e-1}}
=\psi^{p^{e-1}}(a^{p+1})
=\psi^{p^{e-1}}(a)\psi^{p^{e-1}}(a)^p\\
&=\psi^{p^{e-1}}(a)(\psi^p\psi^{p^{e-1}}(a)+p\theta^p\psi^{p^{e-1}}(a))\\
&=\psi^{p^{e-1}}(a)\psi^{p^e}(a)+p \psi^{p^{e-1}}(a)\theta^p\psi^{p^{e-1}}(a)\\
\end{split}
\end{equation*}
which vanishes since $\psi^{p^e}(a)=0$ and $p \psi^{p^{e-1}}(a)=0$
by Proposition \ref{two}.
\end{proof}

\section{The example ring.}

We write $K$ for the localization $\ZZ_{(p)}$ of 
the ring of ordinary integers at  $p$, in which every
prime $r\not=p$ is invertible.
It caries a unique structure of $\lambda$-ring such that
all Adams operations $\psi^q$ coincide with the identity map.

We write $R$ for the polynomial ring $K[x,y]$,
and $\phi\colon R\to R$ for the homomorphism given by
$\phi(x)=x^p-py$, $\phi(y)=y^p$.
\begin{proposition}
For every $f\in R$ one has $\phi(f)=f^p\mod pR$.
\end{proposition}

\begin{proof}
The statement is obviously true if $f\in K$ or $f=x$ or $f=y$.
Moreover if it is true for $f$ and $g$ then it is also true for $f+g$ and $fg$.
\end{proof}

The Proposition tells us that we can define a $\theta^p$ -structure on $R$
by putting $\psi^p(f)=\phi(f)$ and $\theta^p(f)=p^{-1}(f^p-\phi(f))$.
One can extend this to a structure of $\lambda$-ring
by declaring $\psi^r(x)=0$ and $\psi^r(y)=0$ for primes $r$
different from $p$.

\begin{definition}
$F_n\in\ZZ[s,t]$ is defined recursively by:
\begin{equation*}
F_0(s,t)=s,\qquad
F_n(s,t)=F_{n-1}(s,t)^p-pF_{n-1}(t,0)
\end{equation*}
\end{definition}

\begin{lemma}
One has
\begin{equation*}
F_n(s,t)=F_{n-1}(s^p-pt,t^p)\
\end{equation*}
\end{lemma}

\begin{proof}
By induction.
For $n=1$ both sides read $s^p-pt$.
If the statement is true for $n$ then
$F_{n+1}(s,t)=F_n(s,t)^p-pF_n(t,0)
=F_{n-1}(s^p-pt,t^p)-pF_{n-1}(t^p,0)
=F_n(s^p-pt,t^p)$.
\end{proof}

The point of this definition is that
\begin{equation*}
\psi^p F_n(x,y)=F_n(\psi^p(x),\psi^p(y))=F_n(x^p-py,y^p)=F_{n+1}(x,y)
\end{equation*}

\begin{definition}
The ideal $J$ of $R$ is the one generated by
the $p^{e-n}F_n(x,y)$ for $0\leq n\leq e$ and by $y^{p^e}$.
\end{definition}

\begin{proposition}
The ideal $J$ is stable under $\theta^p$.
\end{proposition}

\begin{proof}
From the formulas for $\theta^p$ of a sum and $\theta^p$ of a product
it follows that is sufficient to check that the generators of $J$
are mapped to $J$.
For $n<e$ we have 
\begin{equation*}
\begin{split}
\theta^p(p^{e-n}F_n(x,y))
&=p^{-1}\left((p^{e-n}F_n(x,y))^p-\psi^p(p^{e-n}F_n(x,y))\right)\\
&=p^{(e-n)p-1}F_n(x,y)-p^{e-n-1}F_{n+1}(x,y)\\
\end{split}
\end{equation*}
where the second term is obviously in $J$
and the first term is in $J$ since $(e-n)p-1\geq e-n$.
Furthermore
\begin{equation*}
\begin{split}
\theta^p(F_e(x,y))
&=p^{-1}\left(F_e(x,y)^p-F_e(x^p-py,y^p)\right)\\
&=p^{-1}\left(F_e(x,y)^p-F_{e+1}(x,y))\right)\\
&=p^{-1}(pF_e(y,0))=f_e(y,0)=y^{p^e}
\end{split}
\end{equation*}
and finally $\theta^p(y^{p^e})=0$ since
$\psi^p(y)=y^p$.
 \end{proof}

The Proposition tells us that the quotient ring $A=R/J$
inherits a structure of $\theta^p$-ring,
and in fact a structure of $\lambda$-ring.
We will show now that the class $a$ of $x$ in $R/J$ 
satisfies $a^{p^e+p^{e-1}-1}\not=0$.

The main ingredient for proving this is the following elementary Proposition:
\begin{proposition}
\label{module}
Let $k$ be a commutative ring.
Let $\rho\colon k[\xi,\eta]\to k[x,y]$ be the homomorphism given by
$\rho(\xi)=x^p$, $\rho(\eta)=y^p$.
Let $M$ be an ideal of $k[\xi,\eta]$ with generators $g_j$,
and let $N$ be the ideal of $k[x,y]$ generated bu the $\rho(g_j)$.
Then $k[x,y]/N$ is a free module over $k[\xi,\eta]/M$ via $\rho$
with basis the classes of the $x^ky^m$ for $0\leq k,m<p$.
\end{proposition}

In particular if the class of $\xi^n$ is nonzero in $k[\xi,\eta]/M$
then the class of $x^{pn+p-1}$ is nonzero in $k[x,y]/N$.

\begin{proof}
If $b\in k[x,y]/N$ then it is the class $[f]$ of a certain $f\in k[x,y]$.
Now $f\in k[x,y]$ can be written uniquely as
$\sum_{0\leq k,m<p} x^ky^m \rho(f_{km})$ for certain $f_{km}\in k[\xi,\eta]$,
and therefore $b=\sum_{k,m} [x^ky^m][\rho(f_{km})]$.

On the other hand if  $\sum_{k,m} [x^ky^m][\rho(f_{km})]=0$
then there must be $h_j\in k[x,y]$ such that 
$\sum_{k,m} x^ky^m \rho(f_{km})=\sum h_j \rho(g_j)$.
Moreover each $h_j$ can be written as
$\sum_{k,m} x^ky^m\rho(h_{jkm})$.
Therefore
\begin{equation*}
\sum_{k,m} x^ky^m \rho(f_{km})=\sum_{k,m}x^ky^m\rho(\sum_j g_j h_{jkm})
\end{equation*}
By uniqueness this implies $f_{km}=\sum_j g_j h_{jkm}$.
This means that $[f_{km}]=0$.
\end{proof}

In order to apply this in an induction to prove the main theorem
we need another property of the $F_n$:
\begin{lemma}
\label{power}
If $n>0$ then $F_n(s,t)=F_{n-1}(s^p,t^p)\mod p^n$.
\end{lemma}

\begin{proof}
If $n=1$ this reads $s^p-pt=s^p\mod p$.
If the statement is true for $n$ then
\begin{equation*}
\begin{split}
F_{n+1}(s,t)
&=F_n(s,t)^p-pF_n(t,0)\\
&=(F_{n-1}(s^p,t^p)\mod p^n)^p
-p(F_{n-1}(t^p,0)\mod p^n)\\
&=F_{n-1}(s^p,t^p)^p-p F_{n-1}(t^p,0)\mod p^{n+1}\\
&=F_n(s^p,t^p)
\end{split}
\end{equation*}
Here we used the fact that $u=v\mod p^n$ implies $u^p=v^p\mod p^{n+1}$.
\end{proof}

In order to show that $a^{p^e+p^{e-1}-1}$ is nonzero in
$A=K[x,y]/J$ we prove that is even nonzero in $A/p^{e+1}$
(which is of course not a $\theta^p$-ring):
\begin{theorem}
The class of $x^{p^e+p^{e-1}-1}$ is nonzero
in $\ZZ/p^{e+1}[x,y]$ modulo the ideal generated
by the $p^{e-n}F_n(x,y)$ for $0\leq n\leq e$ and  by $y^{p^e}$
\end{theorem}

\begin{proof}
We use induction in $e$.
For $e=1$ the statement says that $x^p\not=0$ in
$\ZZ/p^2[x,y]/(px,x^p-py,y^p)$.
This is obvious by inspection.

If the statement is true for certain $e$ then by  Proposition \ref{module} 
the class of $x^{p(p^e+p^{e-1}-1)+p-1}=x^{p^{e+1}+p^e-1}$ is nonzero
in $\ZZ/p^{e+1}[x,y]$ modulo the ideal generated
by the $p^{e-n}F_n(x^p,y^p)$ for $0\leq n\leq e$ and  by $y^{p^{e+1}}$.
By Lemma \ref{power} this ideal coincides with the ideal generated
by the $p^{e-n}F_{n+1}(x,y)$ for $0\leq n\leq e$ and  by $y^{p^{e+1}}$.
 
A fortiori the class of $x^{p^{e+1}+p^e-1}$ is nonzero in
$\ZZ/p^{e+2}[x,y]$ modulo the ideal generated by the same 
polynomials and $p^{e+1}F_0(x,y)$.
But this is just the statement for $e+1$.
\end{proof}



\end{document}